\documentclass[12pt, leqno]{article}
\usepackage{amsmath, amsfonts, amsthm}
\begin{document}
\author{Ajai Choudhry}
\title{A diophantine system}
\date{}
\maketitle

\begin{abstract} 
In this paper we find a parametric solution to the hitherto unsolved problem of finding three positive integers such that their sum, the sum of their squares and the sum of their cubes are simultaneously  perfect squares.
\end{abstract}



\bigskip
\bigskip

This paper is concerned with the problem of finding three positive integers such that  the sum of the integers, the sum of the squares of the integers and the sum of the cubes of the integers are all perfect squares. It seems that this  problem has not yet been considered in the literature.

If we denote the three desired  integers by $a,\,b$ and $c$, we need to solve the diophantine system  defined by the following three simultaneous equations:
\begin{align}
a+b+c&=u^2, \label{eq1}\\
a^2+b^2+c^2&=v^2, \label{eq2}\\
a^3+b^3+c^3&=w^2, \label{eq3}
\end{align}
where $u,\,v$ and $w$ are some integers. 

We note that if we take $c=0$, then the left-hand side of Eq.~\eqref{eq3} has a linear factor $a+b$, and it suffices to find rational points on the elliptic curve defined by the pair of quadrics $a^2+b^2=v^2$ and $a^2-ab+b^2=w^2$. This problem has been completely solved by Bremner \cite{Br}. 

We further note that if $a,\,b,\,c$ are any three rational numbers satisfying   the simultaneous equations \eqref{eq1}, \eqref{eq2}, \eqref{eq3}, and  $k$ is any nonzero rational number, then $k^2a,\,k^2b,\,k^2c$ is also a solution of these equations. Thus, any solution in positive rational numbers yields, on appropriate scaling, a solution of these equations in positive integers. It therefore suffices to obtain a  solution of our diophantine system in positive rational numbers.

To solve the above diophantine system, we observe that it follows from \eqref{eq1} and \eqref{eq3} that
\begin{equation}
\begin{aligned}
(a+b+c)^3-(a^3+b^3+c^3)&=u^6-w^2,\\
\mbox{\rm or,}\quad \quad \quad 3(a+b)(b+c)(a+c) &= (u^3-w)(u^3+w).
\end{aligned}
\label{eq13}
\end{equation}
Thus there exist nonzero rational numbers $p$ and $q$ such that 
\begin{equation}
a+b=p(u^3-w),\quad b+c=q(u^3+w),\quad c+a=1/(3pq), \label{eq4}
\end{equation}
and  on solving equations \eqref{eq1} and \eqref{eq4} for $a,\,b,\,c$ and $w$, we get,
\begin{equation}
\begin{aligned}
a&=-(6p^2qu^3-3p(p+q)u^2+1)/(3p(p-q)),\\
b&=(3pqu^2-1)/(3pq),\\
c&=(6pq^2u^3-3q(p+q)u^2+1)/(3q(p-q)),\\
w&=(3pq(p+q)u^3-6pqu^2+1)/(3pq(p-q)).
\end{aligned}
\label{valabc}
\end{equation}

On substituting the values of $a,\,b,\,c$ given by \eqref{valabc} in Eq.~\eqref{eq2} and writing $v=y/\{3pq(p-q)\}$, we get
\begin{multline}
9q^2u^4(8q^2u^2-8qu+3)p^4-6u^2q(12q^3u^3-3q^2u^2-2qu+2)p^3\\
+(27q^4u^4+12q^3u^3+2)p^2-2q(6q^2u^2+1)p+2q^2=y^2. \label{ec1}
\end{multline}

On further substituting $u=1$,  $y=Y/\{4(m^2-8m+8)^2\}$ and
\begin{equation}
q=(m^2-4m-8)/\{2(m^2-8m+8)\}, \label{valq}
\end{equation}
in Eq.~\eqref{ec1}, we get
\begin{multline}
Y^2=36(m^2-8m+24)^2(m^2-4m-8)^2p^4-12(m^2-4m-8)\\ \times (7m^6-136m^5+1112m^4-5120m^3+14272m^2-19968m+1536)p^3\\
+(83m^8-1936m^7+18112m^6-90496m^5+291200m^4\\-705536m^3+1060864m^2-352256m+143360)p^2\\
-8(m^2-8m+8)(m^2-4m-8)(5m^4-56m^3+160m^2-64m+320)p\\
+8(m^2-4m-8)^2(m^2-8m+8)^2.
\label{ec2}
\end{multline}

Now Eq.~\eqref{ec2}, considered as a quartic equation in $p$ and $Y$, represents  a quartic model of an elliptic curve over the field $\mathbb{Q}(m)$. Since the coefficient of $p^4$ in Eq.~\eqref{ec2} is a perfect square,  it is readily found that if we take 
\begin{multline}
p=-(m^{16}-96m^{15}+2688m^{14}-27904m^{13}-108288m^{12}\\
+6494208m^{11}-87138304m^{10}+674709504m^9-3417415680m^8\\
+11595022336m^7-25774522368m^6+35770073088m^5-34652291072m^4\\
+57252249600m^3-130157641728m^2+154014842880m-63870861312)\\
\times (m^2-4m-8)^{-1}(m^2-12m+24)^{-1}(m^2-8m+24)^{-2}\\
\times (m^8-36m^7+672m^6-6944m^5+39936m^4\\
-128256m^3+235520m^2-288768m+258048)^{-1}, \label{valp}
\end{multline}
the right-hand side of Eq.~\eqref{ec2} becomes a perfect square and we thus get a rational point $P$ on the curve \eqref{ec2}. We have omitted the $Y$-coordinate of the point $P$ since it is too cumbersome to write and is also not required for a solution of the problem.

Now, using the relations \eqref{valabc}, we obtain the following values of $a,\,b$ and $c$ satisfying  the simultaneous equations \eqref{eq1}, \eqref{eq2} and \eqref{eq3}:
\begin{equation}
\begin{aligned}
a& = -\{(6q-3)p^2-3pq+1\}/\{3p(p-q)\}, \\
b& = (3pq-1)/(3pq), \\
c& = \{3(2q-1)pq-3q^2+1\}/\{3q(p-q)\}
\end{aligned}
\label{solabc}
\end{equation}
where the values of $p$ and $q$ given by \eqref{valp} and \eqref{valq} respectively. By appropriate scaling, we readily obtain a solution of our  diophantine system in terms of polynomials of degree 68 in the parameter $m$. As these polynomials are too cumbersome to write, we do not give them explicitly.

When $m=3/2$, the above solution yields, on appropriate scaling, the following three integers such that their sum, the sum of their squares and the sum of their cubes are perfect squares:
\[
\begin{aligned}
 &22104703132724392891974197260485203180817980456068478,&\\
& 45051218517398331420875516790921404601474342024364969,&\\
& 273836695120684015976157268469007404280872671207701754. &
\end{aligned}
\]

In fact, it is readily seen that \eqref{solabc} gives a solution in positive rational numbers of the simultaneous equations \eqref{eq1}, \eqref{eq2} and \eqref{eq3} when $m$ is any rational number such that $1.47 < m <1.58$. It immediately follows that we can obtain infinitely many triads of positive integers such their sum, the sum of their squares and the sum of their cubes are perfect squares.

As the solutions obtained by the formulae \eqref{solabc} are large,  we explored the existence of numerically smaller solutions by performing computer trials. Excluding triads of integers that have $k^2$ as a common factor for some integer $k$, we could obtain only two triads of positive integers with the desired property in the range $a+b+c < 10000$. These two triads are as follows:
\[ 
[108, 124, 129] \quad \mbox{\rm and} \quad   [34, 2134, 2873].
\]

Finally we note that the rational point $P$, already found  on the elliptic curve \eqref{ec2} defined over the field $\mathbb{Q}(m)$, is a point of infinite order. We can, therefore, find infinitely many rational points on the elliptic curve \eqref{ec2} and we thus obtain infinitely many parametric solutions of the diophantine system defined by the equations \eqref{eq1}, \eqref{eq2} and \eqref{eq3}. These parametric solutions will yield more integer solutions of our problem.

\begin{center}
\Large
Acknowledgments
\end{center}
 
I wish to  thank the Harish-Chandra Research Institute, Prayagraj for providing me with all necessary facilities that have helped me to pursue my research work in mathematics.

\medskip

\noindent Postal Address: Ajai Choudhry, 
\newline \hspace{1.05 in}
13/4 A Clay Square,
\newline \hspace{1.05 in} Lucknow - 226001, INDIA.
\newline \noindent  E-mail: ajaic203@yahoo.com


\begin{thebibliography}{9}

\bibitem{Br} A. Bremner, A diophantine system, Internat. J. Math. \& Math. Sci. {\bf 9}  (1986), 413--415.

\end{thebibliography}
\end{document}